# Bilateral Solution Bounds and Successive Estimation of Boundedness and Stability Regions for Vector Delay Nonlinear Time-Varying Systems

Mark A. Pinsky

**Abstract.** Stability and boundedness analysis for vector nonlinear systems with variable delays and coefficients remains challenging due to the conservatism of existing methods. Moreover, estimates of the transient behavior of solution norms remain insufficiently developed. This paper presents an approach to estimate the temporal evolution of solution norms and applies it to the analysis of boundedness and stability of vector nonlinear systems with variable delays and coefficients. The method is based on a novel scheme for successive approximations of the original solutions, complemented by the estimates of the corresponding residual norms. This leads to the construction of a scalar nonlinear delay equation whose solutions provide upper bounds for the evolution of residual norms. As a result, bilateral bounds on the original solution norms are obtained, yielding effective boundedness and stability criteria and enabling estimation of the associated regions.

Simulations demonstrate that the proposed approximations rapidly approach the reference boundaries of the regions of interest as the iteration count increases. Moreover, the bilateral bounds progressively approach each other and the norm of the reference solution when the initial function remains within the considered regions.

## 1. Introduction

Vector nonlinear systems with variable time delays and coefficients arise in numerous applications, and analysis of their boundedness and stability is crucial for predicting their temporal behavior.

Over the years, extensive research has addressed some of these problems and provided valuable insights into the temporal evolution of such systems. Numerous monographs [4], [6], [7], [9], [15]–[17], [19]–[22], [23], [26], [28], [30], [31] and countless research papers have examined different aspects of delay system dynamics and developed analytical frameworks for their study. Broader discussions of recent advances and remaining challenges can be found in the review papers [2], [10], [18], [24].

Significant contributions to the stability analysis of delay systems were introduced by N. Krasovskii [29] and B. Razumikhin [43], who extended Lyapunov-functions methods on delayed dynamics. These techniques, based on the construction of Lyapunov-like functionals or specific Lyapunov functions, remain effective analytical tools and have been adapted to certain classes of time-varying delay systems, see, e.g., [2], [6], [7], [32], [49]. However, identifying suitable functionals or functions for general vector delay nonlinear systems remains a major challenge. Nonetheless, for linear delay systems with constant coefficients, these methods provide computationally efficient approaches utilizing linear matrix inequalities (LMI) [16], [18], [19], [27], [45].

Alternative techniques contributing to the stability analysis of delayed systems include, for instance, Halanay's inequality [15], [16], [20], [33] and the Bohl-Perron theorem [4], [5], [17], and other methods. Despite these advances, stability assessments for nonlinear systems with multiple coupled delay equations frequently yield conservative results with limited practical relevance.

Another challenging problem is determining the boundedness of solutions to nonlinear vector delay systems under external perturbations, see, i.e., [8]. This problem has been partly tackled within the framework of input-to-state stability (ISS), which typically relies on restrictive assumptions about the stability of the associated homogeneous systems. For an in-depth discussion on this topic, see the survey [10] and related studies cited herein.

Estimating stability or boundedness regions for delay systems remains another outstanding problem. In principle, these estimates can be obtained using level sets of Lyapunov-Krasovskii functionals or Lyapunov-Razumikhin functions, as discussed in [9], [11], [13], [14], [46]. However, applying these techniques to vector delay nonlinear systems frequently leads to overly conservative and cumbersome criteria with limited practical value. Attempts to simplify these conditions using the norm of the history function [34], [48] often increase conservatism of these techniques.

An effort to estimate stability regions for certain delay nonlinear systems by combining simulations and bifurcation techniques was made in [44]. Fractal structures of stability region boundaries were reported in [1], [11], [12] and [47].

Although boundedness and stability address long-term behavior, estimating transient system dynamics is often of equal practical concern, yet remains insufficiently explored. Moreover, instability phenomena in delay systems have received comparatively limited attention; see, for example, [35], [36].

In [37]–[40], we proposed an approach to stability and boundedness analysis of multidimensional systems of ODEs based on estimating the time evolution of solution norms rather than constructing Lyapunov functions. The method derives scalar auxiliary equations whose solutions upper bounds the corresponding solution norms of the original system. Combined with successive approximations, it yields increasingly accurate estimates of stability and boundedness regions [40]. The approach was extended to nonlinear vector delay systems in [41] and further generalized in [42].

This paper builds on our previous work and develops efficient recursive techniques yielding bilateral solution bounds that are used for estimating boundedness and stability regions of a broad class of vector systems with significant nonlinear components and variable delays, and coefficients. The approach estimates the time evolution of solution norms without relying on Lyapunov-type methodologies. It is based on a novel successive approximation scheme for the original system solutions and on analysis of the temporal evolution of the associated residual terms. We define finite-dimensional approximations of the boundedness and stability regions that enclose trapping or stability balls and introduce novel criteria that yield progressively refined region estimates and bilateral bounds on solution norms.

Numerical simulations show that these estimates converge toward the reference trapping and stability boundaries as the number of iterations increases. Similarly, the bilateral norm bounds approach the reference solution norm remaining within the estimated regions.

The remainder of the paper is organized as follows. Section 2 introduces notation, mathematical preliminaries, and the system formulation. Sections 3 and 4 present the proposed methodology and its theoretical background. Section 5 develops new boundedness and stability criteria. Sections 6 and 7 describe simulation results for delay systems with polynomial and more general nonlinearities. Section 8 concludes the paper and outlines directions for future research.

## 2. Notation, Preliminaries, and Governing Equation

**2.1. Notation.** First, we acknowledge that symbols $\mathbb{R}$, $\mathbb{R}_{\geq 0}$, $\mathbb{R}_+$ and $\mathbb{R}^n$ represent the sets of real, nonnegative and positive real numbers, and $n$-dimensional real vectors; $\mathbb{N}$ is a set of real positive integers, $\mathbb{R}^{n\times n}$ is a set of $n\times n$-matrices, $I \in \mathbb{R}^{n\times n}$ is the identity matrix, $\mathbb{C}$ is a set of complex values, $\mathbb{C}^n$ and $\mathbb{C}^{n\times n}$ are the sets of $n$-dimensional vectors and $n\times n$-matrices with complex entries. Next, we admit that $C([a,b];\mathbb{R}^n)$, $C([a,b]\times\mathbb{R}^n;\mathbb{R}^n)$, $C([a,b];\mathbb{R}_+)$, $C([a,b];\mathbb{R}_{\geq 0})$, and $C([a,b]\times\mathbb{R}^n_{\geq 0};\mathbb{R}_{\geq 0})$ are the spaces of real continuous in all variables functions $\zeta:[a,b]\to\mathbb{R}_+$, $\zeta:[a,b]\times\mathbb{R}^n\to\mathbb{R}^n$, $\zeta:[a,b]\to\mathbb{R}_{\geq 0}$ or $\zeta:[a,b]\times\mathbb{R}^n_{\geq}\to\mathbb{R}_{\geq}$, respectively; the supremum norm $\|\zeta\| \coloneqq \sup_{t\in[a,b]}|\zeta(t)|$, where $|\cdot|$ stands for the Euclidean norm of a vector or the induced norm of a matrix and $b$ can be infinite. Subsequently, we define that the unite step function, $u(t)=1, \forall t \geq 0; u(t)=0, \forall t < 0$ and $\pi(t_i,t_j)=u(t-t_i)-u(t-t_j), t_i < t_j$. Additionally, we acknowledge that that $|x|_\infty = \sup_{i=1,\ldots,n}(|x_i|)$ and $|x|_1 = \sum_{i=1}^n |x_i|, \forall x \in \mathbb{R}^n$; and the upper right-hand derivative in $t$, $D^+ x(t) \coloneqq \dot{x}(t)$.

**2.2.1. Preliminaries.** This paper works with a vector nonlinear system with variable delays and coefficients:
$$\dot{x} = B(t)x + f(t,x(t),x(t-h_1(t)),\ldots,x(t-h_m(t))) + F(t), \forall t \geq t_0 \qquad (2.1)$$
$$x(t) = \varphi(t), \forall t \in [t_0 - \bar{h}, t_0]$$

where $x \in \mathrm{N} \subset \mathbb{R}^n$, $\mathrm{N}$ is a connected subset of $\mathbb{R}^n$, $0 \in \mathrm{N}$, $f(t,\chi_1,\ldots,\chi_{m+1}) \in \mathbb{R}^{n(m+1)+1}$, $\chi_i \in \mathbb{R}^n, i=2,\ldots,m+1$ is a continuous function in all variables that is also locally Lipschitz in $\chi_1$, $f(t,0)=0$, $\varphi \in C([t_0-\bar{h},t_0];\mathbb{R}^n)$, $\bar{h} > 0$, $\|\varphi\| \leq \bar{\varphi}$, matrix $B \in C\left([t_0, \infty); \mathbb{R}^{n \times n}\right)$, $F(t) = F_0 e(t)$, $e \in C\left([t_0, \infty); \mathbb{R}^n\right)$, $\|e(t)\| = 1$, $F_0 \in \mathbb{R}_{\geq 0}$, and scalar functions $h_i \in C\left([t_0, \infty); \mathbb{R}_+\right)$ comply with the following condition:

$$\max_i \sup_{\forall t \geq t_0} h_i(t) = \bar{h} < \infty, \quad \min_i \inf_{\forall t \geq t_0} h_i(t) \geq \underline{h} > 0, \quad i = 1, \ldots, m \tag{2.2}$$

The sequel of this paper assumes that (2.1) admits a unique solution and uses the abridged notation $x(t, \varphi) \coloneqq x(t, t_0, \varphi)$, $\forall t \geq t_0$.

**2.2.2.** Next, we outline some standard definitions of boundedness, stability, and instability of solutions to (2.1) that will be referenced afterwards, see, e.g., [28].

**Definition 1**. Assume that (2.1) admits a unique solution $\forall \varphi \in J \subset C\left(\left[t_0 - \bar{h}, t_0\right]; \mathbb{R}^n\right)$, $\|\varphi\| \leq \bar{\varphi}$. Then, the trivial solution of equation (2.1) with $F_0 \equiv 0$ is called: (1.1) Stable for the set value of $t_0$ if $\forall \varepsilon \in \mathbb{R}_+$, $\exists \delta_1(t_0, \varepsilon) \in \mathbb{R}_+$ such that $|x(t, t_0, \varphi)| < \varepsilon$, $\forall t \geq t_0$ if $\forall \|\varphi\| < \delta_1(t_0, \varepsilon)$. Otherwise, the trivial solution is called unstable. (1.2) Uniformly stable if in the above definition $\delta_1(t_0, \varepsilon) = \delta_2(\varepsilon)$. (1.3) asymptotically stable if it is stable for a given value of $t_0$ and $\exists \delta_3(t_0) \in \mathbb{R}_+$ such that $\lim_{t \to \infty} |x(t, t_0, \varphi)| = 0$ if $\forall \|\varphi\| < \delta_3(t_0)$. (1.4) Uniformly asymptotically stable if it is uniformly stable and in the previous definition $\delta_3(t_0) = \delta_4 = const$.

**Definition 2**. Assume that (2.1) admits a unique solution $\forall \varphi \in J \subset C\left(\left[t_0 - \bar{h}, t_0\right]; \mathbb{R}^n\right)$, $\|\varphi\| \leq \bar{\varphi}$. Then, the trivial solution of equation (2.1) is called: (2.1) Bounded for the set values of $t_0$ if $\exists \delta_6(t_0) \in \mathbb{R}_+$ and $\exists \varepsilon_*(\delta_6) \in \mathbb{R}_+$ such that $|x(t, t_0, \varphi)| < \varepsilon_*$, $\forall t \geq t_0$ if $\|\varphi\| < \delta_6(t_0)$. Otherwise, solution $x(t, t_0, \varphi)$ is called unbounded. (2.2) Uniformly bounded if $\delta_6(t_0)$ is independent upon $t_0$.

**2.2.3.** In [41] we define zero-centered balls $B_{R_i} \coloneqq \varphi \in C\left(\left[t_0 - \bar{h}, t_0\right], \mathbb{R}^n\right) : \|\varphi\| \leq R_i$, $i = 3, \ldots, 6$ with the largest radiuses $R_i$ that contain history functions originating stable or bounded solutions to (2.1), where $R_i = \sup \delta_i$, $i = 3, \ldots, 6$ denote the supremum values of $\delta_i$ for which the corresponding conditions in Definitions 1 or 2 are satisfied.

However, in principle, bounded or stable solutions to (2.1) or its homogeneous counterpart can be generated by $\varphi \notin B_{R_i}$. To account for this, the asymptotic stability region can be defined as a connected subset of the space of continuous history functions such that, $SR \coloneqq \varphi \in C\left(\left[t_0 - s_1, t_0\right], \mathbb{R}^n\right) : \lim_{t \to \infty} |x(t, t_0, \varphi)| \to 0$, where $x(t, t_0, \varphi)$ is a solution of (2.1) with $F_0 = 0$. Similarly, the boundedness region is defined as a connected subset of the space of continuous history functions such that, $BR \coloneqq \varphi \in C\left(\left[t_0 - s_1, t_0\right], \mathbb{R}^n\right):$ $|x(t, t_0, \varphi)| < \varepsilon > 0$, $\forall t \geq t_0$, where $x(t, t_0, \varphi)$ is a solution of (2.1).

It is worth noting that the stability region is frequently termed the region of attraction [28].

Although the above definitions extend naturally to uniform asymptotic stability and boundedness, the numerical estimation of the corresponding regions is substantially more complex.

In principle, numerical simulation of SR and BR requires a finite-dimensional representation of $\varphi(t)$, which can be obtained by projecting the infinite-dimensional space of continuous history functions onto a finite-dimensional subspace using, for example, truncated Fourier series, splines, or related approximation schemes. The resulting systems of nonlinear time-varying ordinary differential equations approximate the original delay differential equations, and their SR and BR boundaries can be estimated using techniques described in [40]. However, applying such methodologies generally involve assessing the correspondence between initial regions and their infinite dimensional approximations. Moreover, improving approximation accuracy within this framework typically requires increased system dimension, leading to a substantial growth in computational complexity.

In contrast, we introduce below an alternative finite-dimensional representation of SR and BR that avoids such dimension growth and allows the original DDE models to be simulated in a form structurally analogous to the ODE-based frameworks developed in our previous work [40].

First, we define a vector $\varphi_s = [\varphi_{s1},...,\varphi_{sn}]^T \in \mathbb{R}^n$, where $\varphi_{sj} = \sup_{t \geq t_0} \varphi_j(t)$, $j = 1,...,n$.

**Definition 3.**

(3.1) The stability domain, $D_S \subset \mathbb{R}^n$ is a connected subset of $\mathbb{R}^n$ that contains the zero vector and all vectors $\varphi_s \in \mathbb{R}^n$ for which $\lim_{t \to \infty} |x(t,t_0,\varphi)| \to 0 : \varphi(t) \equiv \varphi_s$, $\forall t \in [t_0 - \bar{h}, t_0]$, where $x(t,t_0,\varphi)$ is a solution to (2.1) with $F_0 = 0$.

(3.2) The assessed stability region, $ASR := \varphi \in C([t_0 - s_1, t_0], \mathbb{R}^n) : \varphi_s \in D_S$.

(3.3) The boundedness domain, $D_B \subset \mathbb{R}^n$, is a connected subset of $\mathbb{R}^n$ that contains the zero vector and all vectors for which $|x(t,t_0,\varphi)| < \varepsilon$, $\forall t \geq t_0 : \varphi(t) \equiv \varphi_s$, $\forall t \in [t_0 - \bar{h}, t_0]$, where $x(t,t_0,\varphi)$ is a solution to (2.1) with $F_0 > 0$.

(3.4) The assessed boundedness region, $ABR := \varphi \in C([t_0 - s_1, t_0], \mathbb{R}^n) : \varphi_s \in D_B$.

Thus, the definitions (3.1) and (3.3) describe connected regions of $\mathbb{R}^n$ that contain the zero vector and all vectors $\varphi_s \in \mathbb{R}^n$ for which, for given values of $t_0$ and $\varphi(t) \equiv \varphi_s$, $\forall t \in [t_0 - \bar{h}, t_0]$, the solution of (2.1) is either asymptotically stable ($F_0 = 0$) or bounded ($F_0 > 0$). The corresponding sets of $ASR$ and $ABR$ define the sets of continuous history functions that satisfy the respective properties. Ones the boundaries of finite-dimensional regions $D_S$ and $D_b$ are estimated, it becomes straightforward to determine whether a given history function belongs to $ASR / ABR$ or lies outside of them.

We also note that $R_3 \in D_S$ and $R_5 \in D_B$, $B_{R_3} \subset ASR$, and $B_{R_5} \subset ABR$. Thus, Definition 3 can be viewed as a natural generalization of the definitions of stability and boundedness balls and the corresponding radii.

Finally, Definition 3 can be extended to characterize uniform asymptotic stability, exponential stability, or uniform boundedness. However, numerical evaluation of the corresponding regions introduces additional challenges and is therefore not addressed here.

**2.2.4.** Next, we recall some results from our previous work that are used in the sequel. To this end, we introduce a scalar auxiliary equation associated with system (2.1) in the following form [41].

$$\dot{u} = p(t)u + c(t)\left(L(t,u(t),u(t-h_1(t)),...,u(t-h_m(t))) + |F(t)|\right)$$
$$u(t,|\varphi(t)|) = |\varphi(t)|, \forall t \in [t_0 - \bar{h}, t_0] \tag{2.3}$$

where $u(t) \in \mathbb{R}_{\geq 0}$, $p(t) = d(\ln|w(t)|)/dt$, $w(t)$ is the fundamental solution matrix for linear equation $\dot{x} = B(t)x$, $c(t) = |w(t)||w^{-1}(t)|$ is the running condition number of $w(t)$, and a scalar function $L(t,\xi_1,...,\xi_{m+1}) \in \mathbb{R}_{\geq 0}$, $\forall t \geq t_0$, $\xi_i \in \mathbb{R}_{\geq 0}$, $i = 1,...,m+1$ is continuous in all variables, $L(t,0) = 0$ is defined through a nonlinear extension of the Lipschitz continuity condition as follows:

$$|f(t,\chi_1,...,\chi_{m+1})| \leq L(t,|\chi_1|,...,|\chi_{m+1}|), \forall t \geq t_0, \chi_i \in \mathbb{R}^n, i = 1,...,m+1,$$
$$\chi = [\chi_1,...,\chi_{m+1}]^T \in \Omega \subset \mathbb{R}^{n(m+1)} \tag{2.4}$$

where $\Omega$ is a compact and connected subset of $\mathbb{R}^{n(m+1)}$ containing zero. Appendix defines $L$ in close form when $f(t,\chi_1,...,\chi_{m+1})$ is a polynomial, power series or a certain class of rational function in $\chi_1,...,\chi_{m+1}$. Clearly, $L$ is a linear function in $|\chi_i|$ if $f$ is a linear function in $\chi_i$, $i=1,...,m+1$. Note that in all these cases $\Omega \equiv \mathbb{R}^{n(m+1)}$.

Henceforth, for simplicity, we assume that $\Omega \equiv \mathbb{R}^{n(m+1)}$ in the reminder of this paper.

Subsequently, we proved in [41] that $|x(t,\varphi)| \leq u(t,|\varphi|)$, $\forall t \geq t_0$ which enables characterization of the boundedness and stability of the vector equation (2.1) and its homogeneous analog through assessing the corresponding behavior of its scalar counterpart (2.3).

Nonetheless, simulating of $|w(t)|$ over long time intervals can be computationally demanding for large values of $n$ and, in certain cases, $c(t) \to \infty$ as $t \to \infty$ making our estimates overly conservative.

These concerns were addressed in [42] which also develops both the upper and lower bounds for $|x(t,\varphi)|$ under condition that $B(t) = A(t_0) + G(t)$, where $A(t_0) = \limsup\limits_{t \to \infty \, \forall t \geq t_0}(t-t_0)^{-1}\int_{t_0}^{t} B(s)ds \in \mathbb{R}^{n \times n}$ is a nonzero matrix of the general position and $G \in C([t_0,\infty);\mathbb{R}^{n \times n})$ is a zero mean matrix. Let us recall that under this condition $A(t_0)$ assumes simple eigenvalues for all but possibly some isolated values of $t_0 \in \mathbb{R}$ [3]. To proceed, (2.1) was mapped into the eigenbasis of $A$ yielding the following equation:

$$\dot{x}_\times = (A_\times(t_0) + G_\times(t))x_\times + f_\times(t, x_\times(t), x_\times(t-h_1(t)),..., x_\times(t-h_m(t))) + F_\times(t), \forall t \geq t_0$$

$$x_\times(t,\eta) = \eta \equiv V^{-1}\varphi(t), \forall t \in [t_0 - \bar{h}, t_0]$$

where $x_\times = V^{-1}(t_0)x$, $V(t_0) \in \mathbb{C}^{n \times n}$ is the eigenmatrix of $A$ $A_\times = V^{-1}AV = diag(\lambda_1,...,\lambda_n)$, $\lambda_i$, $i=1,...,n$ are the eigenvalues of $A$, $\text{Re}\lambda_i = \alpha_i$, $\alpha_1 > \alpha_i$, $i>1$, $G_\times(t) = V^{-1}G(t)V$, $F_\times = V^{-1}F$ and

$$f_\times(t, x_\times(t), x_\times(t-h_1(t)),..., x_\times(t-h_m(t))) = V^{-1}f(t, Vx_\times(t), Vx_\times(t-h_1(t)),..., Vx_\times(t-h_m(t))).$$

Then, $L_\times$ is defined by applying (2.4) to $f_\times$ as follows: $|f_\times(t, \chi_1,...,\chi_{m+1})| \leq L_\times(t, |\chi_1|,...,|\chi_{m+1}|)$, $\chi_i \in \mathbb{R}^n$, $i=1,...,m+1$, $\chi = [\chi_1,...,\chi_{m+1}]^T \in \Omega \equiv \mathbb{R}^{n(m+1)}$, see a relevant example in Section 6.

Consequently, [42] developed the following bilateral bounds on the original solution norm:

$$\frac{1}{|V^{-1}|}z(t,\phi) \leq |x(t,\varphi(t))| \leq |V|Z(t,\phi), \forall t \geq t_0$$

$$\phi(t) = |V^{-1}\varphi(t)|, \forall t \in [t_0 - \bar{h}, t_0]$$

(2.5)

where $z(t,\phi)$ and $Z(t,\phi)$ are the solutions of the following scalar equations:

$$\dot{Z} = (\alpha_1 + |g(t)|)Z + L_\times(t, Z(t), Z(t-h_1(t)),..., Z(t-h_m(t))) + |F_\times(t)|$$

$$Z(t,\phi) = |V^{-1}\varphi(t)|, \forall t \in [t_0 - \bar{h}, t_0]$$

$$\dot{z} = (\alpha_n - |g(t)|)z - (L_\times(t, z(t), z(t-h_1(t)),..., z(t-h_m(t))) + |F_\times(t)|)$$

$$z(t,\phi) = |V^{-1}\varphi(t)|, \forall t \in [t_0 - \bar{h}, t_0]$$

(2.6)

Here $g(t) = G_\times(t) - iD(t)$, $G_\times(t) = V^{-1}G(t)V$, $D \in \mathbb{R}^{n \times n}$, $D := \text{Im}(diagG_\times)$, and $\alpha_1$ and $\alpha_n$ are the largest and the smallest real parts of the eigenvalues of $A$.

Nonetheless, inequality (2.5) provides relatively accurate estimates for solutions originating from the central portions of the boundedness and stability regions of system (2.1), however, its accuracy deteriorates near their boundaries, leading to conservative assessments of these sets.

These limitations are addressed in the remainder of this paper by representing the solution of the original equation as the sum of an approximate solution and the corresponding residual. The norm of the residual is then estimated using the right-hand side of (2.5), yielding bilateral bounds that progressively approach one another for solutions whose histories lie within the boundedness or stability regions.

## 3. Successive Approximations

To be consistent with simulations provided in Section 6 and 7, equation (2.1) is rewritten as

$$\dot{x} = \left(A(t_0) + G(t)\right)x + E(t)x(t - h_0(t)) + f\left(t, x(t), x(t - h_1(t)), \ldots, x(t - h_m(t))\right) + F(t), \ \forall t \geq t_0$$
$$x(t, \varphi) = \varphi(t), \ \forall t \in \left[t_0 - \bar{h}, t_0\right] \tag{3.1}$$

where function $f$ now contains no linear terms, matrices $A(t_0)$, $G(t)$, $E \in C\left([t_0, \infty); \mathbb{R}^{n \times n}\right)$, functions $f, F, h_i, i = 1, \ldots, m$ are defined above, matrix and $h_0 \in C\left([t_0, \infty); \mathbb{R}_+\right)$ complies with (2.2). Next, we assume that equation (3.1) assumes a unique solution $\forall \|\varphi\| \leq \bar{\phi}$ which is represented as follows:

$$x(t, \varphi) = Y_K(t, \varphi) + z_K(t, \varphi), \ t \geq t_0, \ 1 \leq K \in \mathbb{N} \tag{3.2}$$

where an approximate solution: $Y_m := [t_0, \infty) \times J \to \mathbb{R}^n$, $Y_m(t, \varphi) = \varphi(t)$, $\forall t \in \left[t_0 - \bar{h}, t_0\right]$, $J \subset C\left(\left[t_0 - \bar{h}, t_0\right]; \mathbb{R}^n\right)$, and $z_K := [t_0, \infty) \times J \to \mathbb{R}^n$, $z_K(t, \varphi) = 0$, $\forall t \in \left[t_0 - \bar{h}, t_0\right]$ is the residual of this approximation. It will be clear subsequently that the right side of the differential equation for $z_K$ includes $\varphi(t)$ which justifies our prior notation, i.e., $z_K = z_K(t, \varphi)$.

Then, applying the triangular inequality and its reverse part to (3.2) yields the following bilateral bounds on $|x(t, \varphi)|$:

$$|Y_K(t, \varphi)| - |z_K(t, \varphi)| \leq |x(t, \varphi)| \leq |Y_K(t, \varphi)| + |z_K(t, \varphi)|, \ t \geq t_0 \tag{3.3}$$

where the left side of (3.3) should be adjusted to zero if it takes negative values.

In turn, we represent the approximate solution in the following form:

$$Y_K(t, \varphi) = \sum_{k=1}^{K} y_k(t, \varphi), \ 1 \leq K \in \mathbb{N} \tag{3.4}$$

where $y_i := [t_0, \infty) \times J \to \mathbb{R}^n$ solve the following linear equations:

$$\dot{y}_1 = A(t_0)y_1 + F(t), \ \forall t \geq t_0$$
$$y_1(t, \varphi) = \varphi(t), \ \forall t \in \left[t_0 - \bar{h}, t_0\right]$$
$$\dot{y}_2 = A(t_0)y_2(t) + G(t)y_1(t) + E(t)y_1(t - h_0(t)) + f\left(t, y_1(t), y_1(t - h_1(t)), \ldots, y_1(t - h_m(t))\right), \ \forall t \geq t_0$$
$$y_2(t, \varphi) = 0, \ \forall t \in \left[t_0 - \bar{h}, t_0\right] \tag{3.5}$$
$$\dot{y}_k = A(t_0)y_k(t) + G(t)y_{k-1}(t) + E(t)y_{k-1}(t - h_0(t)) + f\left(t, Y_{k-1}(t), Y_{k-1}(t - h_1(t)), \ldots, Y_{k-1}(t - h_m(t))\right) -$$
$$f\left(t, Y_{k-2}(t), Y_{k-2}(t - h_1(t)), \ldots, Y_{k-2}(t - h_m(t))\right), \ \forall t \geq t_0$$
$$y_k(t, \varphi) = 0, \ \forall t \in \left[t_0 - \bar{h}, t_0\right], \ 3 \leq k \leq K$$

where, for brevity, we adopt that $Y_k(t, \varphi) := Y_k(t), \ 2 < k \leq K$.

System (3.5) comprises a sequence of linear ODEs with constant coefficients that are recursively coupled through time-varying delay terms. Since the coupling terms are explicitly determined by solutions of the preceding equations, the system admits a nearly closed-form representation of solutions, facilitating both qualitative analysis and numerical evaluation.

In turn, the residual term for (3.5) is defined through the following nonlinear delay equation:

$$\dot{z}_K = Az_K + G(t)z_K + E(t)z_K(t-h_0(t)) + G(t)y_K(t) + E(t)y_K(t-h_0(t)) +$$
$$f(t, z_K(t)+Y_K(t), z_K(t-h_1(t))+Y_K(t-h_1(t)),...,z_K(t-h_m(t))+Y_K(t-h_m(t))) -$$
$$f(t, Y_{K-1}(t), Y_{K-1}(t-h_1(t)),...,Y_{K-1}(t-h_m(t)))$$
$$z_K(t,\varphi) = 0, \forall t \in [t_0 - \bar{h}, t_0], \forall K \geq 1$$

where $f(t, Y_0(t), Y_0(t-h_1(t)),...,Y_0(t-h_m(t))) \equiv 0$.

We assume that the latter equations admit unique solutions. Instead of solving them directly, we estimate the corresponding solution norm using the right-hand side of inequality (2.5). Following the procedure used in the derivation of (2.5), we transform this equation into the eigenbasis of $A$, which yields the following equation:

$$\dot{\hat{z}}_K = A_\times \hat{z}_K + G_\times(t)\hat{z}_K + E_\times(t)\hat{z}_K(t-h_0(t)) + \Gamma_K(t, \hat{z}_K(t),...,\hat{z}_K(t-h_1(t)),...,\hat{z}_K(t-h_m(t)))$$
$$\hat{z}_K(t,\varphi) = 0, \forall t \in [t_0 - \bar{h}, t_0], \forall K \geq 1 \quad (3.6)$$

where $\hat{z}_K = V^{-1}z_K$, $E_\times(t) = V^{-1}E(t)V$,

$$\Gamma_K(t, \hat{z}_K(t),...,\hat{z}_K(t-h_1(t)),...,\hat{z}_K(t-h_m(t))) =$$
$$V^{-1}\begin{pmatrix} \hat{f}(t, \hat{z}_K(t)+Y_K(t), \hat{z}_K(t-h_1(t))+Y_K(t-h_1(t)),...,\hat{z}_K(t-h_m(t))+Y_K(t-h_m(t))) - \\ f(t, Y_{K-1}(t), Y_{K-1}(t-h_1(t)),...,Y_{K-1}(t-h_m(t))) + G(t)y_K(t) + E(t)y_K(t-h_0(t)) \end{pmatrix}$$

and

$$\hat{f}(t, \hat{z}_K(t)+Y_K(t), \hat{z}_K(t-h_1(t))+Y_K(t-h_1(t)),...,\hat{z}_K(t-h_m(t))+Y_K(t-h_m(t))) =$$
$$f(t, V\hat{z}_K(t)+Y_K(t), V\hat{z}_K(t-h_1(t))+Y_K(t-h_1(t)),...,V\hat{z}_K(t-h_m(t))+Y_K(t-h_m(t)))$$

In turn, applying (2.4) to $\Gamma_K(t, \chi_1,...,\chi_{m+1})$, yields the following inequality:

$$|\Gamma_K(t, \chi_1,...,\chi_{m+1})| \leq L_{\Gamma_K}(t, |\chi_1|,...,|\chi_{m+1}|), \forall \chi = [\chi_1,...,\chi_{m+1}]^T \in \mathbb{R}^{n(m+1)} \quad (3.7)$$

where $L_{\Gamma_K} \in C([t_0, \infty) \times \mathbb{R}_{\geq 0}^{m+1}; \mathbb{R}_{\geq 0})$, $\forall K \geq 1$. Note that the next section shows the application of this inequality to some systems with polynomial nonlinearities.

To upper bound $|\hat{z}_K(t,\varphi(t))|$, we set up a scalar counterpart for a vector equation (3.6) as follows:

$$\dot{Z}_K = (\alpha_1 + |g(t)|)Z_K + |E_\times(t)|Z_K(t-h_0(t)) + L_{\Gamma_K}(t, Z_K(t), Z_K(t-h_1(t)),...,Z_K(t-h_m(t)))$$
$$Z_K(t,\varphi) = 0, \forall t \in [t_0 - \bar{h}, t_0], \forall K \geq 1 \quad (3.8)$$

Note that the right-side of (3.8) implicitly depends upon $\varphi(t)$ since $Y_k = Y_k(t,\varphi), k \geq 1$. Next, we assume that equation (3.8) admits a unique solution $\forall t \geq t_0, \forall \|\varphi\| \leq \bar{\phi}, \forall K \geq 1$. Then, an application of the upper bound in (2.5) yields the following inequality:

$$|z_K(t,\varphi(t))| \leq |V|Z_K(t,\varphi), \forall t \geq t_0; K \geq 1$$
$$Z_K(t,\varphi) = 0, \forall t \in [t_0 - \bar{h}, t_0] \quad (3.9)$$

where $Z_K(t,\varphi)$ is a solution to a scalar equation (3.8).

Lastly, we write (3.3) in the following form:

$$|Y_K(t,\varphi)| - |V|Z_K(t,\varphi) \leq |x(t,\varphi)| \leq |Y_K(t,\varphi)| + |V|Z_K(t,\varphi), t \geq t_0, \forall \|\varphi\| \leq \bar{\phi}, \forall K \geq 1 \quad (3.10)$$

## 4. Behavior of Approximate Solutions Over Long Time Intervals

The key results of this study are based on the proposition stated below.

**Theorem 1.** Assume that $f \in C\left([t_0,\infty) \times \mathbb{R}^{n(m+1)}; \mathbb{R}^n\right)$, $f(t,0)=0$, and a local Lipschitz-like condition holds in all variables of this function starting from the second, i.e.:

$$\left|f(t,\chi_1,...,\chi_{m+1}) - f(t,\chi_1^*,...,\chi_{m+1}^*)\right| \leq \sum_{i=1}^{m+1} q_i(t,\tilde{\chi})|\chi_i - \chi_i^*|, \ \forall |\chi - \chi^*| \leq \tilde{\chi} > 0, \ \forall t \geq t_0, \ \chi_i, \chi_i^* \in \mathbb{R}^n,$$

$$\chi, \chi^* = \left[\chi_1^*,...,\chi_{m+1}^*\right]^T \in \mathbb{R}^{n(m+1)}, \ q_i \in C\left([t_0,\infty) \times \mathbb{R}_{\geq 0}; \mathbb{R}_{\geq 0}\right), \ q_i(t,\tilde{\chi}) \leq q(\tilde{\chi}) \in \mathbb{R}_+, \ \forall t \geq t_0, \ i=1,...,m \quad (4.1)$$

Also, we assume that $\varphi \in C\left([t_0-\bar{h}, t_0]; \mathbb{R}^n\right)$, $\bar{h} > 0$, $\|\varphi\| \leq \bar{\varphi} \in \mathbb{R}_{\geq 0}$, $F_*(t) = F_0 e(t)$, $e \in C\left([t_0,\infty); \mathbb{R}^n\right)$, $\|e(t)\| = 1$, $F_0 \in \mathbb{R}_{\geq 0}$, $h_i \in C\left([t_0,\infty); \mathbb{R}_+\right)$, $i=0,...,m+1$ comply with (2.2), matrices $G, E \in C\left([t_0,\infty); \mathbb{R}^{n \times n}\right)$, $|G(t)|, |E(t)| < \infty$, $A(t_0) \in \mathbb{R}^{n \times n}$ is a Hurwitz matrix of general position.

Then,
$$|Y_K(t)| \leq (\|\varphi\| + O(F_0)) P_{K-1}(t) \exp(-\eta t) + O(F_0), \forall t \geq t_0, K > 1 \quad (4.2)$$

where $P_{K-1}(t)$ is a polynomial of degree $K-1$, and $\alpha_1$ is the largest real part of the eigenvalue of $A$, $-\eta = \alpha_1$.

**Proof.** First, we recall that $|\exp(A(t-\tau))| \leq \exp(-\eta(t-\tau))$. Next, all equations in (3.5) admit a unique solution given by the variation-of-parameters. Standard norm inequalities applied to the corresponding solution of the first equation of this sequence yield:

$$|y_1(t,\varphi)| \leq \left[\|\varphi\| e^{-\eta(t-t_0)} + (F_0/v_1)(1 - e^{-\eta(t-t_0)})\right] = O(\|\varphi\| + F_0), \ \forall t \geq t_0$$
$$|y_1(t,\varphi)| = |\varphi(t)| \leq \|\varphi\|, \ \forall t \in \left[t_0 - \bar{h}, t_0\right] \quad (4.3)$$

Consequently, using (3.3), (4.1) and (4.3), we can write that:

$$|y_2(t,\varphi)| \leq e^{-\eta t} \int_{t_0}^{t} e^{\eta \tau}\left(\|G\| |y_1(\tau)| + \|E\| y_1(\tau - h_0(\tau)) + q(y_{1,s})\left(|y_1(\tau)| + \sum_{i=1}^{m+1} |y_1(\tau - h_i(\tau))|\right)\right) d\tau.$$ Next, we use that,

$|y_1(\tau - h_i(\tau))| \leq \|\varphi\|$, $\forall t \in \left[t_0 - \bar{h}, t_0\right]$, $i = 0,...,m+1$ and combining this with the estimate:

$\exp(h_{i,inf}) \leq \exp(h_i(t)) \leq \exp(h_{i,s})$, $\forall t \geq t_0$, $i = 0,...,m+1$, where $h_{i,s} = \sup h_i(t)$, $\forall t \geq t_0$ and $h_{i,Inf} = \inf h_i(t)$, $\forall t \geq t_0$. Applying standard integration rules then yields,

$|y_2(t)| \leq (\|\varphi\| + O(F_0)) P_1(t) \exp(-\eta t) + O(F_0)$, $\forall t \geq t_0$, $|y_2(t)| = |\varphi(t)|$, $\forall t \in \left[t_0 - \bar{h}, t_0\right]$, where $P_1(t)$ is a first degree polynomial in $t$.

In order, we assume by induction that:

$|y_K(t)| \leq (\|\varphi\| + O(F_0)) P_{K-1}(t) \exp(-\eta t) + O(F_0)$, $\forall t \geq t_0$, $|y_K(t)| = 0$, $\forall t \in \left[t_0 - \bar{h}, t_0\right]$, $K > 2$, where $P_{K-1}(t)$ is a polynomial of degree $K-1$. Applying the same reasoning, shows that the induction hypothesis extends to $K := K+1$ □

Consequently, (4.2) and (4.3) imply that

$$\lim_{t \to \infty} |Y_K(t,\varphi)| = O(F_0), \ \forall K \geq 1 \quad (4.4)$$

Lastly, combining (4.4) with the right-hand side of (3.10) yields that:

$$\lim_{t \to \infty} |x(t,\varphi)| \leq |V| \lim_{t \to \infty} Z_K(t,\varphi), \ F_0 = 0, \ K \geq 1$$
$$\lim_{t \to \infty} |x(t,\varphi)| \leq O(F_0) + |V| \lim_{t \to \infty} Z_K(t,\varphi), \ F_0 > 0, \ K \geq 1 \quad (4.5)$$

## 5. Boundedness and Stability Criteria

Thus, the boundedness and stability analysis of the original vector equation (3.1) is reduced to the study of the scalar equation (3.8), as reflected in the following statements.

**Theorem 2**. Assume that conditions of Theorem 1 hold, $F_0 = 0$, and $\exists K \geq 1$ for which a scalar function $L_{\Gamma_K} \in C\left([t_0, \infty) \times \mathbb{R}_\geq^n; \mathbb{R}_\geq\right)$ is defined by (3.7). Suppose that for this value of $K$ equation (3.8) admits a unique solution $\forall \|\varphi\| \leq \bar{\phi}$, and the trivial solution of this equation is either stable for the set value of $t_0$, uniformly stable, asymptotically stable for the set value of $t_0$, or uniformly asymptotically stable.

Then, the trivial solution of (3.1) satisfies the corresponding stability properties. Furthermore, the asymptotic stability region of equation (3.1) contains the asymptotic stability region for the corresponding equation (3.8).

**Proof**. The proof readily follows from the first part of inequality (4.5) □

**Theorem 3**. Assume that conditions of Theorem 1 hold, $F_0 > 0$, and $\exists K \geq 1$ for which a scalar function $L_{\Gamma_K} \in C\left([t_0, \infty) \times \mathbb{R}_\geq^n; \mathbb{R}_\geq\right)$ is defined by (3.7). Suppose that that for this value of $K$ equation (3.8) admits a unique solution $\forall \|\varphi\| \leq \bar{\phi}$ and the solutions of the scalar equation (3.8) are bounded for the set value of $t_0$ or uniformly bounded provided that the history functions, $\varphi$, belongs to the corresponding boundedness regions of this equation.

Then, the solutions of the corresponding vector equation (3.1) are bounded for the same value of $t_0$ or uniformly bounded, respectively. Moreover, boundedness regions for equation (3.1) contain the corresponding boundedness regions for equation (3.8).

**Proof**. The proof follows immediately from the second part of inequality (4.5) □

Inequalities (3.3) yield bilateral bounds on the original solution norm, which converge to the true solution as the number of iterations increases, provided that the initial function belongs to the boundedness or stability regions. Their computation relies on nearly closed-form estimates delivered by equations (3.3) and on the simulation of the scalar equation (3.8) and is therefore essentially independent of the dimension of the original system.

While the long-time evolution of $|Y_K(t,\varphi)|$ is governed by inequality (4.2), its magnitude increases rapidly as the initial function approaches the boundaries of the boundedness and stability domains from within. This transient behavior enables efficient estimation of these regions while bypassing direct simulation of equation (3.8).

## 6. Simulations

This section applies the developed approach to estimate the boundaries of boundedness and stability domains ($D_S$ and $D_b$) for two delay systems with time-dependent, nonperiodic coefficients and dissipative or conservative nonlinearities, enabling direct assessment of whether a given history function belongs to the assessed boundedness or stability regions ($ASR / ABR$); see Definition 3. Each system consists of two coupled oscillators resembling Van der Pol–like and Duffing–like dynamics with variable coefficients, external perturbations, and delays. These models arise in various applications (see, e.g., [25]).

Although sufficient conditions for local boundedness and stability can be derived using methods developed in our recent works [41], [42] and in the related studies cited above, these approaches typically yield conservative estimates, particularly for trapping regions. Moreover, the techniques proposed in [41] are not applicable here since, for both considered systems, $\lim_{t \to \infty} c(t) = \infty$.

### 6.1 Scalar Equations for Systems of Van der–Pol–like and Duffing-like Oscillators with Delays

Let us assume that in (3.1)

$$A = \begin{pmatrix} 0 & 1 & 0 & 0 \\ -(\omega_1^2 + d) & -\alpha_1 & d & 0 \\ 0 & 0 & 0 & 1 \\ d & 0 & -(\omega_2^2 + d) & -\alpha_2 \end{pmatrix}, G = \begin{pmatrix} 0 & 0 & 0 & 0 \\ -g_{21}(t) & 0 & 0 & 0 \\ 0 & 0 & 0 & 0 \\ 0 & 0 & -g_{43}(t) & 0 \end{pmatrix}, E(t) = G(t)$$

$$F(t) = \begin{pmatrix} 0 & F_0 \sin \omega_0 t & 0 & 0 \end{pmatrix}^T$$

where $\alpha_1 = 0.4$, $\alpha_2 = 0.2$, $\omega_1^2 = 1$, $\omega_2^2 = 4$, $g_{21} = a_1 \sin r_1 t + a_2 \sin r_2 t$, $g_{43} = b_1 \sin s_1 t + b_2 \sin s_2 t$,

$r_1 = 3.14$, $r_2 = 6.15$, $s_1 = 3.1$, $s_2 = 6.28$, $a_1 = a_2 = b_1 = b_2 = 0.1$, $\omega_0 = 5.43$. Let us also recall that $V = [v_{ij}]$, $i, j = 1, \ldots, 4$ is the eigenvector matrix of $A$.

Additionally, we assume that for Van der Pol-like system:

$$f(t-h_1) = \begin{pmatrix} 0 & -\mu_1 x_2^3(t-h_1) & 0 & -\mu_2 x_4^3(t-h_1) \end{pmatrix}^T$$

and for a Duffing-like system

$$f(t-h_1) = \begin{pmatrix} 0 & -\mu_1 x_1^3(t-h_1) & 0 & -\mu_2 x_3^3(t-h_1) \end{pmatrix}^T$$

Note that the values of $\mu_1$ and $\mu_2$ are indicated on the figures' captions.

Next, we derive a scalar equation (3.5) for two systems indicated above. To this end, we define the following scalars, $\kappa_i = \sum_{k=1}^{4} abs(v_{ik})$, $i = 1, \ldots, 4$, where $abs(\varsigma) = |\varsigma|$, $\forall \varsigma \in \mathbb{R}$, $abs(\varsigma) = (\varsigma_1^2 + \varsigma_2^2)^{1/2}$, $\forall \varsigma = \varsigma_1 + i\varsigma_2 \in \mathbb{C}$.

Subsequently, we observe the following inequalities:

$$abs\left(\sum_{l=1}^{4}(v_{il}\hat{z}_{lK})\right)^j \leq \left(\sum_{l=1}^{4} abs(v_{il}\hat{z}_{lK})\right)^j = \left(\sum_{l=1}^{4} abs(v_{il})abs(\hat{z}_{lK})\right)^j \leq \left(\sum_{l=1}^{4} abs(v_{il}) Z_K\right)^j =$$

$$\left(\sum_{l=1}^{4} abs(v_{il})\right)^j Z_K^j = \kappa_i^j Z_K^j, \quad i = 1, \ldots, 4, \ j = 1, 2, 3, \ K \in \mathbb{N}$$

where we used that $abs(\hat{z}_{lK}) \leq \left\|[z_{1K}, \ldots, z_{4K}]^T\right\| \equiv Z_K$.

Next, we estimate some components in $|\Gamma_K|$ as follows:

$$abc\left(\left(\sum_{l=1}^{4}(v_{il}\hat{z}_{lK}) + Y_{iK}\right)^3 - Y_{i,K-1}^3\right) \leq$$

$$abc\left(\sum_{l=1}^{4}(v_{il}\hat{z}_{lK})\right)^3 + 3abc\left(\sum_{l=1}^{4}(v_{il}\hat{z}_{lK})\right)^2 abc(Y_{iK}) + 3abc\left(\sum_{k=1}^{4}(v_{il}\hat{z}_{lK})\right)Y_{iK}^2 + abc(Y_{iK}^3 - Y_{i,K-1}^3) \leq$$

$$\kappa_i^3 Z_K^3 + 3\kappa_i^2 Z_K^2 \bar{Y}_{iK} + 3\kappa_i Z_K \bar{Y}_{iK}^2 + abc(Y_{iK}^3 - Y_{i,K-1}^3), \ i = 1, \ldots, 4$$

where we assumed that $\bar{Y}_{iK} = abs(Y_{iK})$.

Next, we establish the following inequality:

$$\mu_1 abc\left(\left(\sum_{k=1}^{4}(v_{il}\hat{z}_{lK}) + Y_{iK}\right)^3 - Y_{i,K-1}^3\right) + \mu_2 abc\left(\left(\sum_{k=1}^{4}(v_{jl}\hat{z}_{lK}) + Y_{jK}\right)^3 - Y_{j,K-1}^3\right) \leq$$

$$\mu_1\left(\kappa_i^3 Z_K^3 + 3\kappa_i^2 Z_K^2 \bar{Y}_{iK} + 3\kappa_i Z_K \bar{Y}_{iK}^2 + abc(Y_{iK}^3 - Y_{i,K-1}^3)\right) +$$

$$\mu_2\left(\kappa_j^3 Z_K^3 + 3\kappa_j^2 Z_K^2 \bar{Y}_{jK} + 3\kappa_j Z_K \bar{Y}_{jK}^2 + abc(Y_{jK}^3 - Y_{j,K-1}^3)\right)$$

where it is assumed that if $i = 1$, then $j = 3$ (Van der-Pol-like system), and if $i = 2$, then $j = 4$ (Duffing-like system).

Next, we estimate that,

$$|\Gamma_K| = |V^{-1}(f(t, V\hat{z}_K + Y_K) - f(t, Y_{K-1}))| \leq$$

$$|V^{-1}| \begin{pmatrix} Z_K^3(\mu_1 \kappa_i^3 + \mu_2 \kappa_j^3) + 3Z_K^2(\mu_1 \kappa_i^2 \bar{Y}_{iK} + \mu_2 \kappa_j^2 \bar{Y}_{jK}) + 3Z_K(\mu_1 \kappa_i Y_{iK}^2 + \mu_2 \kappa_j Y_{jK}^2) + \\ + \mu_1 abc(Y_{iK}^3 - Y_{i,K-1}^3) + \mu_2 abc(Y_{jK}^3 - Y_{j,K-1}^3) \end{pmatrix}$$

where, as above, it is assumed that if $i = 1$, then $j = 3$ (Van der-Pol-like system), and if $i = 2$, then $j = 4$ (Duffing-like system).

Lastly, for our systems the scalar equation (3.8) takes the following form:

$$\dot{Z}_1 = \alpha_1 Z_1 + \left(\left|V^{-1}G(t)V\right|\right)\left(Z_1(t) + Z_1(t-h_0)\right) +$$

$$\left|V^{-1}\right| \begin{pmatrix} \left(\mu_1\kappa_2^3 + \varepsilon_2\kappa_4^3\right)Z_1^3(t-h_1) + \\ 3\left(\mu_1\kappa_2^2 abs(y_{2,1}(t-h_1)) + \mu_2\kappa_4^2 abs(y_{4,1}(t-h_1))\right)Z_1^2(t-h_1) + \\ 3\left(\mu_1\kappa_2 y_{2,1}^2(t-h_1) + \mu_2\kappa_4 y_{4,1}^2(t-h_1)\right)Z_1(t-h_1) + \\ \mu_1 abs(y_{2,1}^3(t-h_1)) + \mu_2 abc(y_{4,1}^3(t-h_1)) \end{pmatrix} +$$

$$\left|V^{-1}\left(G(t)(y_1(t) + y_1(t-h_0))\right)\right|$$

$$Z_1(t,\phi) = 0, \forall t \in \left[t_0 - \bar{h}, t_0\right]$$

and

$$\dot{Z}_K = \alpha_1 Z_K(t) + \left|V^{-1}G(t)V\right|\left(Z_K(t) + Z_K(t-h_0)\right) +$$

$$\left|V^{-1}\right| \begin{pmatrix} \left(\mu_1\kappa_2^3 + \varepsilon_2\kappa_4^3\right)Z_K^3(t-h_1) + 3\left(\mu_1\kappa_2^2 \bar{Y}_{2K}(t-h_1) + \mu_2\kappa_4^2 \bar{Y}_{4K}(t-h_1)\right)Z_K^2(t-h_1) + \\ 3\left(\mu_1\kappa_2 Y_{2K}^2(t-h_1) + \mu_2\kappa_4 Y_{4K}^2(t-h_1)\right)Z_K(t-h_1) + \\ \mu_1 abc\left(Y_{2K}^3(t-h_1) - Y_{2,K-1}^3(t-h_1)\right) + \mu_2 abc\left(Y_{4K}^3(t-h_1) - Y_{4,K-1}^3(t-h_1)\right) \end{pmatrix} +$$

$$\left|V^{-1}\left(G(t)(y_K(t) + y_K(t-h_0))\right)\right|,$$

$$Z_K(t,\phi) = 0, \forall t \in \left[t_0 - \bar{h}, t_0\right], \forall K \geq 1$$

(6.1)

**6.2. Simulation of Boundedness and Stability Domains and Bilateral Solution Bounds**

Our simulations rely on two components: (i) Theorems 2–3 with inequalities (3.10) and (4.5), and (ii) closed-form approximate solutions of (3.1) derived from (3.5). In both cases, the simulated estimates approach the boundaries of trapping and stability domains (Definition 3) for fixed T as $K$ increases, whereas the reference boundaries are obtained from independent simulations of (3.1).

The first strategy yields conservative inner approximations, whereas the second typically produces outer approximations, often enabling bilateral bounds on the estimated boundaries. The latter, more efficient, approach relay on simulating nearly closed-form solutions of the linear equations (3.5), while the former - simulates $Z_K(t,\varphi_s), t \in [t_0, T]$ using the scalar equation (6.1) combined with nearly closed-form solutions to linear equations (3.5) with $y_1(t,\varphi_s) = \varphi_s$, $y_i(t,0) = 0, \forall t \in \left[t_0 - \bar{h}, t_0\right]$, $1 < i \leq K$. In both cases, the computational cost is essentially independent of system dimension.

Moreover, inequality (3.10) provides bilateral bounds on the solution norm that converge to the exact norm as iteration counts increases, provided the solution remains within a trapping or stability region.

Although boundedness and Lyapunov stability (Definitions 1–2) are defined on infinite time horizon $t \in [t_0, \infty)$, numerical assessment of system dynamics over a finite interval $t \in [t_0, T]$ is $T$-dependent and often increasingly conservative as $T$ grows.

Simulations used MATLAB's variable-step DDE solver with default tolerances $10^{-5}/10^{-6}$; smaller tolerances (e.g., $10^{-8}$) improved plausibility of the results in some cases.

In the simulations, the components of the initial vector $\varphi_s \in \mathbb{R}^4$ were represented in double polar coordinates as follows:

$$\varphi_{s1}(t_0) = r_1(t_0)\cos(\theta_1(t_0)),\ \varphi_{s2}(t_0) = r_1(t_0)\sin(\theta_1(t_0)),$$
$$\varphi_{s3}(t_0) = r_2(t_0)\cos(\theta_2(t_0)),\ \varphi_{s4}(t_0) = r_2(t_0)\sin(\theta_2(t_0))$$

Both angular coordinates were discretized with step size $\pi/48$. For each fixed angular pair, the radial components were adjusted to approximate region boundaries. Simulations started at small radial values located inside the trapping or stability domains and were progressively increased via a binary search until a prescribed threshold $\varpi$ was exceeded. Increasing/decreasing $\varpi$ typically led to more/ fewer conservative estimates.

The same procedure was used to estimate boundary components of $\varphi_s$ associated with boundedness and stability domains, based on simulations of either the right-hand side of inequality (3.3) or the scalar equation (3.8). Both approaches produced indistinguishable results; accordingly, lower boundary estimates are reported from simulations of (3.8).

Moreover, when the vector $\varphi_s$ is near the reference boundary, the magnitude of $|Y_K(t,\varphi_s)|$ rapidly increases and approaches the same threshold $\varpi$, providing a complementary boundary estimate.

For clarity, estimated boundaries are projected onto the $\varphi_{s1} \times \varphi_{s2}$-plane, as projections onto other coordinate planes yield qualitatively similar results. These projections are shown in the figures in solid, dashed, and dot-dashed lines, corresponding to simulations of the original system (3.1), equation (6.1), and the condition $|Y_K(t,\varphi_s)| = \varpi$, respectively. The plots in this section show that approximation accuracy remains intact for relatively large delays and nonlinear scaling parameters, when the evaluation interval sufficiently exceeds the maximal lag.

**6.3.1. Stability Domains and Solution Bounds for Coupled Delay Van der-Pol-like System**

Figure 1 shows the boundaries of estimated stability-domains near the origin for a coupled delay Van der Pol–like system, computed using four, six, and eight iterations over different time intervals. Except for Figure. 1a, the stability-domain obtained from simulations of the original system (solid line) lies between two estimates: one from (6.2) (dashed line) and another from the time histories of $|Y_m(t,\varphi_s)|$ (dot-dashed line). Figures 1a–1c indicate that, with four iterations and fixed parameters, accuracy varies as $T$ increases. For the same parameter values and $T = 40\,\text{or}\,80$, the accuracy is recovered with six iterations (Figures. 1d–1e) and further improved with eight iterations (Figure 1d), even for larger nonlinear scaling parameters, $\mu_1 = \mu_2 = 60$.

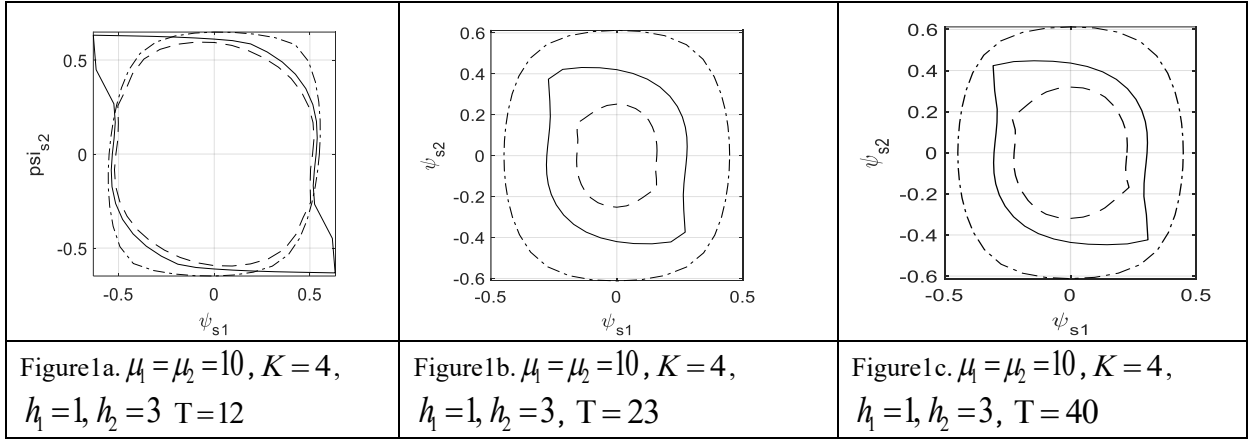

| Figure1a. $\mu_1 = \mu_2 = 10$, $K = 4$, $h_1 = 1, h_2 = 3$ T=12 | Figure1b. $\mu_1 = \mu_2 = 10$, $K = 4$, $h_1 = 1, h_2 = 3$, T=23 | Figure1c. $\mu_1 = \mu_2 = 10$, $K = 4$, $h_1 = 1, h_2 = 3$, T=40 |
|---|---|---|

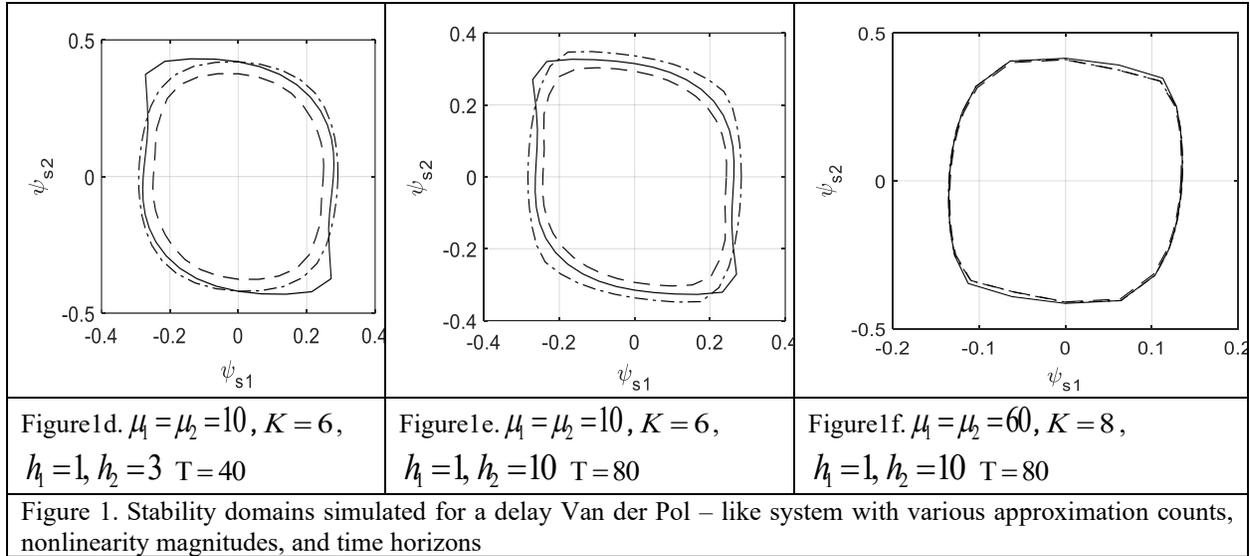

| Figure1d. $\mu_1=\mu_2=10$, $K=6$, $h_1=1, h_2=3$ T=40 | Figure1e. $\mu_1=\mu_2=10$, $K=6$, $h_1=1, h_2=10$ T=80 | Figure1f. $\mu_1=\mu_2=60$, $K=8$, $h_1=1, h_2=10$ T=80 |

Figure 1. Stability domains simulated for a delay Van der Pol – like system with various approximation counts, nonlinearity magnitudes, and time horizons

Figure 2 depicts temporal behavior of upper and lower bonds enclosing the solution norms for the original system of coupled Van-der-Pol-like oscillators (solid line), obtained by applying inequality (3.10). All estimates use six iterations and originate from the same initial vector.

In Figures 2a–2c, the initial vector, $\varphi_s$, lies withing and near the boundary of stability-domain, whereas in Figure 2d, $\varphi_s$ lies outside the estimated stability-domain but within the referenced domain. . In all cases, the reference trajectory remains enclosed by the computed bounds. Moreover, as $K$ increases, the bounds approach each other and the reference curve in Figures 2a-2c.

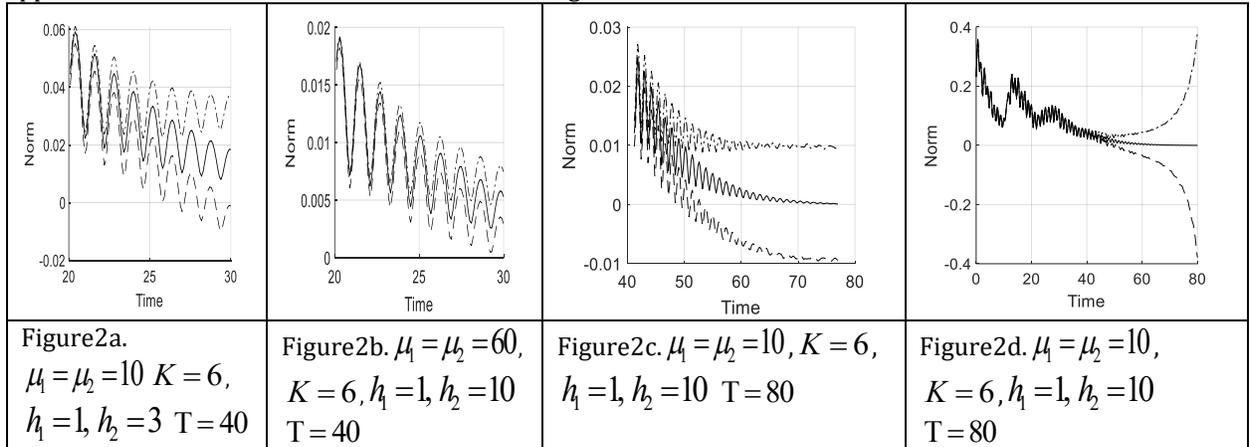

| Figure2a. $\mu_1=\mu_2=10$ $K=6$, $h_1=1, h_2=3$ T=40 | Figure2b. $\mu_1=\mu_2=60$, $K=6, h_1=1, h_2=10$ T=40 | Figure2c. $\mu_1=\mu_2=10$, $K=6$, $h_1=1, h_2=10$ T=80 | Figure2d. $\mu_1=\mu_2=10$, $K=6, h_1=1, h_2=10$ T=80 |

Figure 2. Bilateral solution bounds simulated for a delay Van der Pol–like system with various delays, nonlinearity magnitudes, and time horizons

### 6.3.2. Stability Domains and Solution Bounds for Coupled Delay Duffing -like System

Figure 3 shows estimated stability domains for a system of coupled delay Duffing-like oscillators obtained using six and eight iterations; curve identification follows the conventions of Figures 1 and 2. The first three plots demonstrate acceptable accuracy obtained in six iterations, with further improvement observed at eight iterations and marked refinement of the estimated domain achieved at six iterations for larger values of T=80.

Figure 3d illustrates a representative temporal evolution of the solution norm together with its upper and lower bounds, computed using six iterations. All curves originate from the same initial vector, and the bounds converge toward the reference trajectory as $K$ increases.

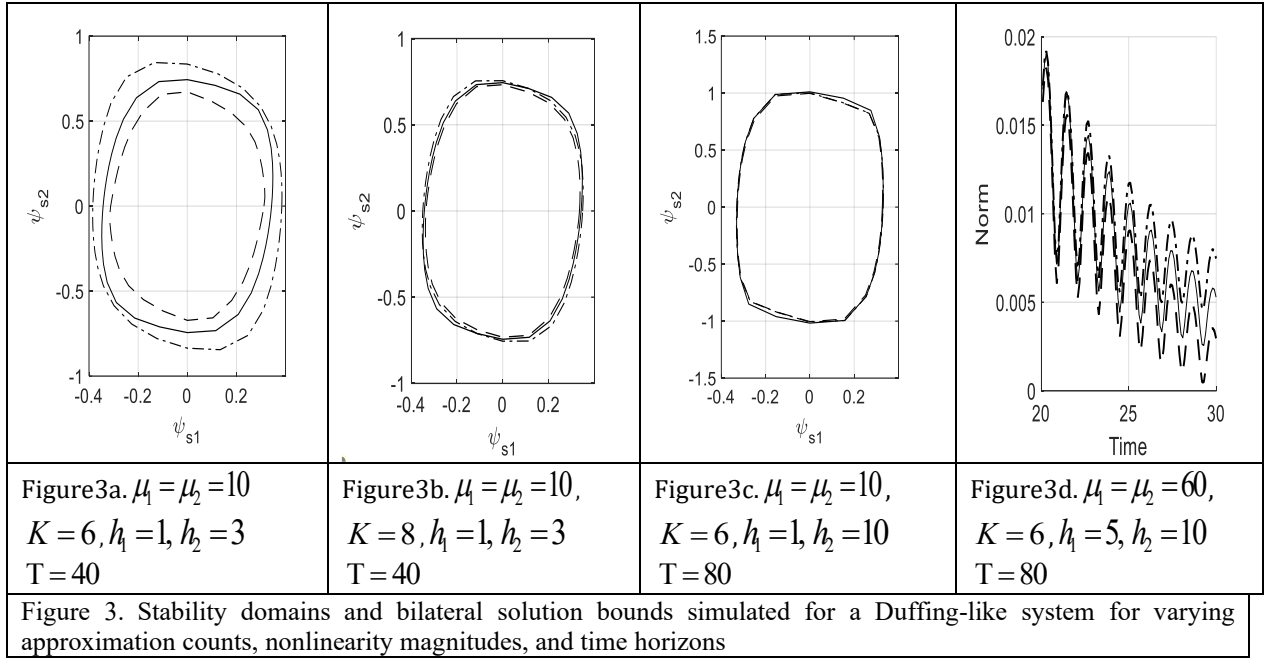

| Figure3a. $\mu_1 = \mu_2 = 10$ | Figure3b. $\mu_1 = \mu_2 = 10$, | Figure3c. $\mu_1 = \mu_2 = 10$, | Figure3d. $\mu_1 = \mu_2 = 60$, |
|---|---|---|---|
| $K = 6$, $h_1 = 1$, $h_2 = 3$ | $K = 8$, $h_1 = 1$, $h_2 = 3$ | $K = 6$, $h_1 = 1$, $h_2 = 10$ | $K = 6$, $h_1 = 5$, $h_2 = 10$ |
| T = 40 | T = 40 | T = 80 | T = 80 |

Figure 3. Stability domains and bilateral solution bounds simulated for a Duffing-like system for varying approximation counts, nonlinearity magnitudes, and time horizons

### 6.3.3. Trapping Domains and Solution Bounds for Nonhomogeneous Van der-Pol-like and Duffing-like Systems with Delay

To our knowledge, the problem of estimating the temporal evolution of solution norms for nonhomogeneous vector nonlinear systems with variable delay and coefficients has received limited attention in the literature. In [41] and [42], we derived local boundedness criteria leading to conservative estimates of boundedness regions. This section demonstrates that the present methodology addresses this gap by effectively capturing the influence of nonhomogeneous terms on system behavior; see Fig. 4, where curve identification follows the conventions of the preceding figures.

Figures 4a and 4b depict the projections of the estimated trapping regions onto the $\varphi_{1,s} \times \varphi_{2,s}$ coordinate plane over eight iterations for nonhomogeneous Van der Pol–like and Duffing–like systems with delays, respectively. In both plots, the dashed and dot-dashed lines closely follow the solid reference curve, with minimal deviation. Subsequent simulations confirm that the estimation errors remain sufficiently small over eight iterations for a broad range of system parameters, provided the simulation interval is chosen to be sufficiently large.

Figures 4c and 4d show the temporal evolution of the reference solution norm and its bilateral bounds for Van der Pol–like and Duffing–like delay systems, respectively, obtained using the proposed methodology. In both cases, the initial vectors lie within the reference boundedness- domain but outside its inner estimate, causing the bounds to diverge while still enclosing the reference solution.

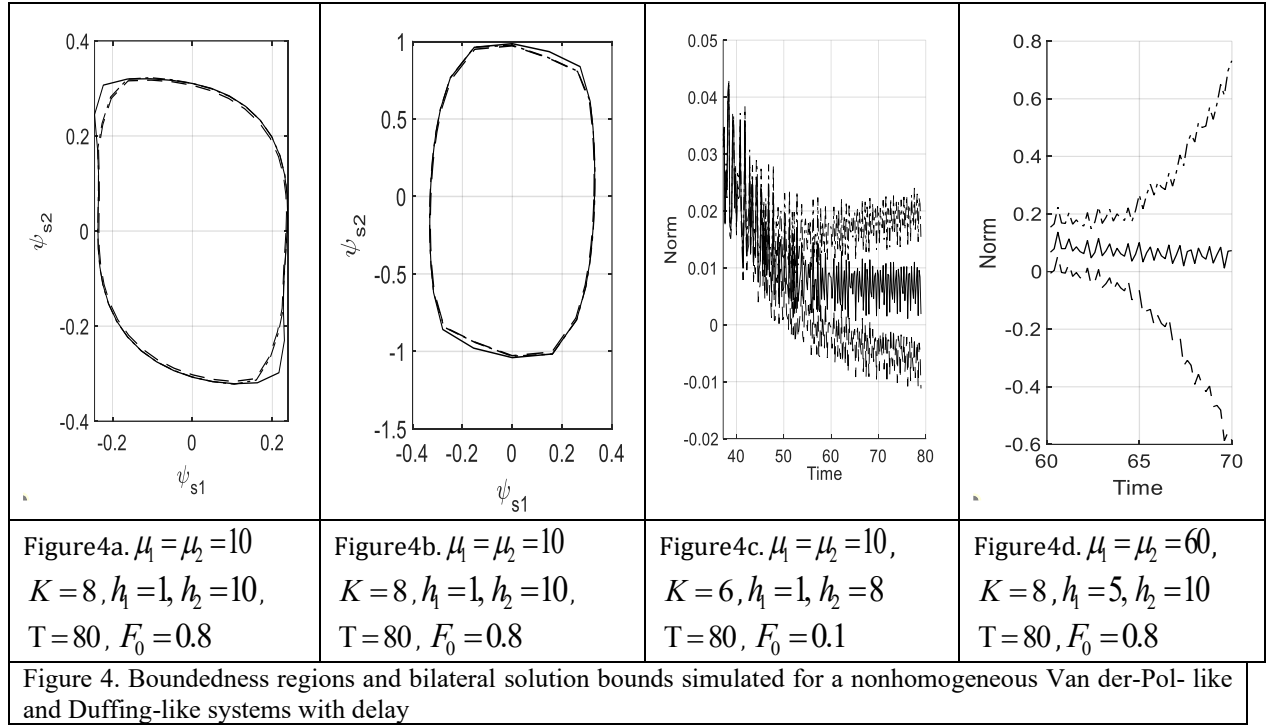

| Figure 4a. $\mu_1 = \mu_2 = 10$, $K=8, h_1=1, h_2=10$, $T=80, F_0=0.8$ | Figure 4b. $\mu_1 = \mu_2 = 10$, $K=8, h_1=1, h_2=10$, $T=80, F_0=0.8$ | Figure 4c. $\mu_1 = \mu_2 = 10$, $K=6, h_1=1, h_2=8$, $T=80, F_0=0.1$ | Figure 4d. $\mu_1 = \mu_2 = 60$, $K=8, h_1=5, h_2=10$, $T=80, F_0=0.8$ |
|---|---|---|---|

Figure 4. Boundedness regions and bilateral solution bounds simulated for a nonhomogeneous Van der-Pol- like and Duffing-like systems with delay

## 7. Estimation of Boundedness and Stability Domains for Delay Systems with Complex Nonlinearities

As noted earlier, the scalar equation (3.8) can be readily derived for vector delay systems whose nonlinear components are polynomials, power series, or certain rational functions. More general nonlinearities may be treated by approximations with such admissible functions; otherwise, deriving (3.8) becomes considerably more involved.

In contrast, the linear equations (3.5) admit nearly closed-form solutions under broad assumptions that render their components as known functions of time. This property enables estimation of boundedness and stability domains for (3.1), even in the presence of complex nonlinear components. For polynomial nonlinearities, this approach was demonstrated in Section 6.3 alongside the primary technique yielding inner approximations of the regions of interest. The complementary strategy relies on simulating components of the initial vector, $\varphi_s$, which induces rapid growth of $\|Y_K(t,\varphi_s)\|$ to a prescribed threshold.

Representative simulations support this strategy for estimating boundedness and stability domains of vector delay systems with complex nonlinearities and variable coefficients. They indicate that the approximation error decreases as $K$ increases, provided $T$ is sufficiently large. While this approach does not ensure inner or outer approximations of the boundedness or stability domains, it offers a valuable starting point for subsequent, more computationally intensive analyses.

This section evaluates the effectiveness of the proposed technique by applying it to two vector delay systems of coupled Van der Pol–like and Duffing-like oscillators augmented with Gaussian and hyperbolic tangent terms, which model impulsive and saturation nonlinearities, respectively. In the first case, the nonlinear components of both systems are represented as follows:

$$f = \begin{pmatrix} 0 & -(\mu_1 x_2^3(t-h_1) + \mu_3 G(\Delta_2)) & 0 & -(\mu_2 x_4^3(t-h_1) + \mu_4 G(\Delta_4)) \end{pmatrix}^T$$
$$f = \begin{pmatrix} 0 & -(\mu_1 x_1^3(t-h_1) + \mu_3 G(\Delta_1)) & 0 & -(\mu_2 x_3^3(t-h_1) + \mu_4 G(\Delta_3)) \end{pmatrix}^T$$

(7.1)

where $G(\Delta) = \dfrac{1}{q\sqrt{2\pi}} \exp\left(-\dfrac{1}{2}\left(\dfrac{\Delta}{q}\right)^2\right)$, $\Delta_1 = x_1(t-h_1) - x_1^0$, $\Delta_2 = x_2(t-h_1) - x_2^0$, $\Delta_3 = x_3(t-h_1) - x_3^0$, $\Delta_4 = x_4(t-h_1) - x_4^0$, $\mu_1 = \mu_2 = -0.1$, $\mu_3 = \mu_4 = -1$, $x_1^0 = x_2^0 = x_3^0 = x_4^0 = 7$ and $h_0 = 1, h_2 = 10, q = 1, T = 80$.

In the latter case, the nonlinear components of the systems are specified as follows:

$$f = \begin{pmatrix} 0 & -\left(\mu_1 x_2^3(t-h_1) + \mu_3 \tanh(k_1 x_2(t-h_1))\right) & 0 & -\left(\mu_2 x_4^3(t-h_1) + \mu_4 \tanh(k_2 x_4(t-h_1))\right) \end{pmatrix}^T$$
$$f = \begin{pmatrix} 0 & -\left(\mu_1 x_1^3(t-h_1) + \mu_3 \tanh(k_1 x_1(t-h_1))\right) & 0 & -\left(\mu_2 x_3^3(t-h_1) + \mu_4 \tanh(k_2 x_3(t-h_1))\right) \end{pmatrix}^T$$
(7.2)

where $\mu_1 = \mu_2 = -0.1$, $\mu_3 = \mu_4 = -1$ and $h_0 = 1, h_1 = 10, k_1 = k_2 = 1, T = 80$.

Note that, in the first case, $f(0) \neq 0$. Consequently, the system does not admit an equilibrium at the origin, and the proposed technique estimates a boundedness domain about the origin of coordinates.

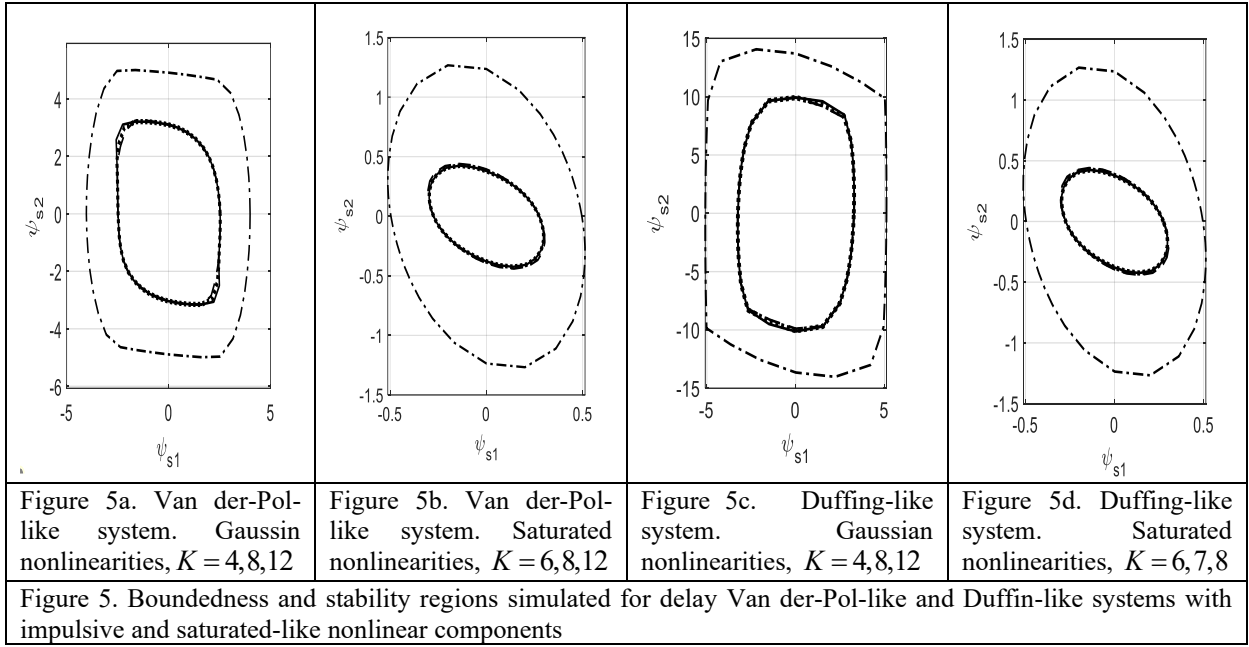

| Figure 5a. Van der-Pol-like system. Gaussin nonlinearities, $K = 4,8,12$ | Figure 5b. Van der-Pol-like system. Saturated nonlinearities, $K = 6,8,12$ | Figure 5c. Duffing-like system. Gaussian nonlinearities, $K = 4,8,12$ | Figure 5d. Duffing-like system. Saturated nonlinearities, $K = 6,7,8$ |
|---|---|---|---|

Figure 5. Boundedness and stability regions simulated for delay Van der-Pol-like and Duffin-like systems with impulsive and saturated-like nonlinear components

Figures 5 show projections onto the $\varphi_1 \times \varphi_2$ plane of estimated boundedness and stability domain boundaries about the origin of coordinates of the coupled systems with nonlinearities given by (7.1) or (7.2). The solid curve denotes the reference boundary obtained from simulations of the original system (3.1), while the dot–dashed, dashed, and bold–dotted curves correspond to simulations of $|Y_K(t,\varphi_s)| = \varpi$ for successive values of $K$ listed in the figure captions.

In all panels, the boundaries obtained for the two largest values of $K$ nearly coincide with each other and with the reference curve, indicating that the proposed methodology provides efficient estimates of boundedness and stability domains for vector delay systems with complex nonlinearities. Appropriate values of $K$ may therefore be selected by comparing estimates obtained for successive values of this parameter.

## 8. Conclusion and Further Research

Estimating sets of initial functions that yield bounded or stable solutions for vector delay nonlinear time-varying systems is a challenging yet practically important problem. Most existing methods provide local stability and boundedness criteria that result in overly complex and conservative conditions with limited practical relavance.

This work extends the methodology developed in [40] for nonautonomous ODEs to vector nonlinear DDEs with variable delays and coefficients. The method is based on recursively approximating the original solutions by solutions of carefully constructed systems of linear ordinary differential equations with constant coefficients and forcing terms, which admit nearly closed-form solutions and preserve the boundedness and stability properties of the original system. The norm of the approximation residual is then upper bounded by the solution of a derived scalar delay equation that implicitly incorporates the initial vector function. This approach replaces direct analysis of the original vector system with analysis of its scalar counterpart, substantially reducing computational complexity regardless of the number of equations coupled in the original DDE.

Within this framework, bilateral bounds for the norms of the original solutions and new boundedness and stability criteria are obtained, yielding inner estimates that converge to the reference boundaries of the boundedness and stability domains as the iteration count increases. Numerical simulations demonstrate rapid error reduction over a broad parameter range and show that the bilateral bounds approach each other and the reference solution norm when the initial function belongs to the boundedness or stability regions.

The simulations further reveal that the linear ordinary differential equations governing successive approximations exhibit rapid growth near the boundaries of the boundedness and stability domains. This transient behavior enables efficient estimation of these regions and provides effective initial approximations for more rigorous, albeit computationally more intensive, methods. The approach remains effective for systems with complex nonlinearities, including impulsive-like and saturation-like components, with simulations confirming rapid convergence to the reference boundaries.

Finally, the methodology extends naturally to nonautonomous vector delay systems with impulsive and stochastic components and can be used to shape boundedness and stability regions in control applications, thereby enhancing robustness assessment.

**Acknowledgment**. The program used to simulate the models in Section 6 was coded by Steve Koblik.

**Appendix**. Let $x = [x_1 \ x_2]^T$ and $f$ is defined, e.g., as follows: $f = [a_1(t) x_1 x_2^2(t-h_1) \quad a_2(t) x_1^3(t-h_2)]^T$, $a_i(t) \in \mathbb{R}$, $i = 1, 2$. Then, $|f|_2 \leq |f|_1 \leq |a_1||x_1||x_2^2(t-h_1)| + |a_2||x_1^3(t-h_2)|$
$\leq |a_1||x(t)||x(t-h_1)|^2 + |a_2||x(t-h_2)|^3$, where we use that $|x_i^n| \leq |x|^n$, $|x_i^n(t-h)| \leq |x(t-h)|^n$, $i = 1, 2, n \in \mathbb{N}$. Clearly, such inequalities can be extended on power series and some rational functions. Really, let $f(x) = (x_1(t-h_1))/(x_1(t-h_2) + x_2(t-h_2))$, $x \in \mathbb{R}^2$, $h_1 \neq h_2 \in \mathbb{R}_+$; then, $|f(x)| \leq |x_1(t-h_1)|/(|x_1(t-h_1)| - |x_2(t-h_2)|)$. The left side of the former inequality is defined and continuous in $t$ if $|x_1(t-h_1)| \neq |x_2(t-h_2)|$, $\forall t \geq t_0$.